\documentclass[12pt,a4paper]{article}
[1994/12/01]
\usepackage{amsmath,amsthm}

\usepackage{amsfonts}
\usepackage{graphics}
\usepackage{graphicx}
\usepackage{amssymb}
\usepackage{amsthm}
\usepackage{pstricks,pst-node,pst-tree,pstricks-add}
\usepackage{newlfont}
\usepackage{epsfig}
\usepackage[final]{pdfpages}

\newtheoremstyle{theorem}
  {10pt}          
  {10pt}  
  {\sl}  
  {\parindent}     
  {\bf}  
  {. }    
  { }    
  {}     
\theoremstyle{theorem}
\newtheorem{theorem}{Theorem}[section]

\newtheoremstyle{defi}
  {10pt}          
  {10pt}  
  {\rm}  
  {\parindent}     
  {\bf}  
  {. }    
  { }    
  {}     
\theoremstyle{defi}

\usepackage{latexsym}
\usepackage{color}
\usepackage{amssymb}
\usepackage{amsmath}
\usepackage{amsthm} \usepackage{graphicx}
\newtheorem{lem}{Lemma}[section] \newtheorem{coro}{Corollary}

\newtheorem{defi}{Definition}[section]
\newtheorem{case}{Case}\newtheorem{rem}{Remark}[section]
\newtheorem{exa}{Example}[section]
 \newtheorem{cla}{Claim}
\newtheorem{ill}{Illustration}[section]

\newcommand{\bl}{\begin{lem}}
\newcommand{\el}{\end{lem}}
\newcommand{\bt}{\begin{theorem}}
\newcommand{\et}{\end{theorem}}
\newcommand{\bc}{\begin{coro}}
\newcommand{\ec}{\end{coro}}
\newcommand{\bd}{\begin{defi}}
\newcommand{\ed}{\end{defi}}
\newcommand{\bp}{\begin{proof}}
\newcommand{\ep}{\end{proof}}
\newcommand{\br}{\begin{rem}}
\newcommand{\er}{\end{rem}}
\newcommand{\be}{\begin{exa}}
\newcommand{\ee}{\end{exa}}
\newcommand{\bca}{\begin{case}}
\newcommand{\eca}{\end{case}}
\newcommand{\bcl}{\begin{cla}}
\newcommand{\ecl}{\end{cla}}
\newcommand{\bil}{\begin{ill}}
\newcommand{\eil}{\end{ill}}

\parindent 18pt\date{}
\begin{document}

\title{Korovkin results and Frobenius optimal
approximants for infinite dimensional bounded
linear operators}
\author{Kiran Kumar, M.N.N.Namboodiri
\thanks{Department of di Mathematics,
``CUSAT'' - Cochin (INDIA) (E-mail: {\tt mnnadri@gmail.com})},
Stefano Serra-Capizzano \thanks{Department 'Fisica e Matematica',
Via Valleggio 11, 22100 Como (ITALY) (E-mail: {\tt stefano.serrac@uninsubria.it})}}
\maketitle
\begin{abstract}
The classical as well as non commutative Korovkin-type theorems deal with convergence of positive linear maps with respect to modes of convergences such as norm convergence and weak operator convergence. In this article, Korovkin-type theorems are proved for convergence of completely positive maps with respect to weak, strong and uniform clustering of sequences  of matrices of growing order. Such modes of convergence were originally considered for Toeplitz matrices (see \cite{st1},\cite{tyr}). As an application,
we translate the Korovkin-type approach used in the finite dimensional case,
in the setting of preconditioning large linear systems with
Toeplitz structure, into the infinite dimensional context of operators
acting on separable Hilbert spaces. The asymptotic of these pre-conditioners are obtained
 and analyzed using the concept of completely positive maps. It is observed that any two limit points of the same sequence of pre-conditioners are the same modulo compact operators. Finally, we prove the generalized versions of the Korovkin type
 theorems in \cite{st1}.
\end{abstract}
\textbf{MSC Classification (1999)} 41A36, 47A58, 47B35
\section{Introduction}
The classical approximation theorem due to Korovkin in \cite{Kor}, unified many approximation processes such as Bernstein polynomial approximation of continuous real functions. Inspired by this discovery, several mathematicians extended the Korovkin's theorems in many ways and to several settings, including function spaces, abstract Banach lattices, Banach algebras, Banach spaces and so on. Such developments are referred to as Korovkin-type approximation theory 
(see \cite{Alt1},\cite{Alt2} and references there reported). The non commutative versions of the Korovkin's theorem, can be found in many papers (see \cite{Alt1},\cite{Ln1},\cite{Ln2},\cite{Ln3},\cite{Uchyama} and references therein and refer to \cite{mnn1}, for new perspectives). In most of these developments, the underlying modes of convergence has been the norm,
strong or weak operator convergence of linear operators in a Hilbert space.

 In this paper, we prove such non commutative Korovkin-type
 theorems with the modes  of convergence induced by strong, weak or uniform clustering of sequences of matrices of growing order. Such notions have  already been used for the special case of Toeplitz matrices in connection with the Frobenius optimal approximation of matrices of large size which has been widely considered in the numerical linear algebra literature for the design of efficient solvers of complicate linear systems of large size (see \cite{st1},\cite{tyr}). More specifically, the approximation is constrained in spaces of low complexity: as examples of high interest in several important applications (see \cite{Chan}, \cite{szego} and references therein), we may mention algebras of matrices associated to fast transforms like Fourier, Trigonometric, Hartley, Wavelet transforms (\cite{Kal}, \cite{CV}) or we may mention spaces with prescribed patterns of sparsity. In the context of general linear systems, accompanied with the minimization in Frobenius norm, these techniques were originally considered and studied by Huckle (see \cite{bto} and references there reported), while the specific adaptation in the Toeplitz context started with the work of Tony Chan in \cite{Tony}. 
 
 More recently, a unified structural analysis was introduced by the third author in connection with the Korovkin theory, which represents a nice branch of the theory of functional approximation. More precisely, the analysis of clustering of the preconditioned systems which gives a measure for
the approximation quality is reduced to classical Korovkin test on a finite number of very elementary symbols associated to equally elementary Toeplitz matrices (Jordan matrices). Here we consider the same approach in an operator theory setting. We translate in the infinite dimensional context of operators acting on separable Hilbert spaces, the Korovkin type approach used in the finite dimensional case in the setting of preconditioning large linear systems with Toeplitz structure arising in various important applications.

 The paper is organized as follows.
In section 2, we give basic definitions and theorems, including the notion of complete positivity, and the description of the problem. In section 3, we prove the new versions of non commutative Korovkin-type theorems. As an application of these general theorems, we consider the example of Frobenius optimal maps and prove some results on their limit points in the next section. Analysis of the LPO sequences in \cite{st1}, is generalized in the fifth section. Finally, a concluding section ends the paper.

\section{Basic definitions and description of the problem}
We begin with the classical Korovkin's theorem.
\bt
Let $\{\Phi_n\}$ be a sequence of positive linear maps on $C[0,1]$.
If
\begin{equation}\nonumber
\Phi_n(f)\rightarrow f \,\,\textrm{for every f in the set}\,\,\{1,x,x^2\},
\end{equation}
then
\begin{equation}\nonumber
\Phi_n(f)\rightarrow f \,\,\textrm{for every f in }\,\,C[0,1].
\end{equation}
\et
Here the convergence is the uniform convergence of sequence of functions.
For the non commutative versions of this theorem, we need the notion of \textit{complete positivity.}
\bd
Let $\mathbb{A}$ and $\mathbb{B}$ be $C^{*}$ algebras with identities $1_{\mathbb{A}}$ and $1_\mathbb{B}$ respectively
and $\Phi :\mathbb{A}\to \mathbb{B}$ be a positive linear map such that $\Phi(1_\mathbb{A})\leq(1_\mathbb{B})$. For each
positive integer n, let $\Phi_n : M_n(\mathbb{A})\to M_n(\mathbb{B})$ be defined as $\Phi_n(a_{i,j})
=(\Phi(a_{i,j} ))$ for every matrix $(a_{i,j})\in M_n(A)$. If
$\Phi_n$ is positive for each n,
then $\Phi$ is called \textit{completely positive}.
\ed
\br
Let $CP\left( {\mathbb{A},{\mathbb{B}}} \right)$ denote the class of all completely positive maps
(CP-maps) $\Phi$ from $\mathbb{A}$ to $\mathbb{B}$ such that $\Phi(1_{\mathbb{A}}) \leq (1_{\mathbb{B}})$.
Then it is well known that $CP\left( {\mathbb{A},\mathbb{B}} \right)$ is compact and convex in the Kadison's
B.W topology\cite{Arv1}.
\er
\br
Let $\mathbb{A},\mathbb{B}$ and $\mathbb{C}$ be three $C^{*}$ algebras. If $\Phi \in CP\left( {\mathbb{A},\mathbb{B}} \right)$
and  $\Psi \in CP\left( {\mathbb{B},\mathbb{C}} \right)$, then 	the composition 
$\Psi\circ\Phi \in CP\left( {\mathbb{A},\mathbb{C}} \right).$ 
\er
We state one of the most fundamental result,
the Stinespring dilation Theorem \cite{stn}.
\bt \label{Theorem 2}
Let $\mathbb{A}$ and $\mathbb{B}$ be $C^{*}$ algebras with identities $1_{\mathbb{A}}$ and $1_{\mathbb{B}}$ respectively.
Let $\Phi :\mathbb{A} \to \mathbb{B}$ be a completely positive linear map such that
$\Phi(1_{\mathbb{A}}) \leq (1_{\mathbb{B}})$. Assume that $\mathbb{B}$
is a sub algebra of $\mathbb{B}(\mathbb{H})$ for some Hilbert space $\mathbb{H}$.
Then there exists a representation $\pi$ of $\mathbb{A}$ on a Hilbert space $\mathbb{K}$ and a
bounded linear map V from $\mathbb{H}$ to $\mathbb{K}$ such that
$ \Phi \left( a \right) = V^{*} \pi \left( a \right)V $ for every $a \in \mathbb{A}.$
\et
\br 
It is known that if either $\mathbb{A}$ or $\mathbb{B}$ is commutative then every positive linear map
is completely positive.
\er
\bd \label{Schwarz map}
A positive linear map $\Phi$ from a $C^{*}$ algebra $\mathbb{A}$ to a $C^{*}$ algebra $\mathbb{B}$ is called a \textit{Schwarz map} if $\Phi(a^*a)\geq \Phi(a)\Phi(a^*)$ for all $a$ in $\mathbb{A}$.
\ed
\br
It can be easily seen that every completely positive map of norm less than 1, is a \textit{Schwarz map}.
\er
The non commutative versions of the classical Korovkin's theorem have been obtained by various researchers for positive maps, Schwarz maps and CP-maps, in the settings $C^*$-algebras and $W^*$-algebras. A short survey of these developments can be found in \cite{mnn1}.
For example, the following is such a version proved in \cite{Ln3}. 

\bt
Let the sequence of Schwarz maps $\Phi_n$, from $\mathbb{A}$ to $\mathbb{B}$ be such that
$\Phi_n(1_{\mathbb{A}})\leq 1_{\mathbb{B}}$. Then the set
\begin{equation}\nonumber
C=\{a \in \mathbb{A};\Phi_n(a)\rightarrow a,\Phi_n(a^*a+aa^*)\rightarrow (a^*a+aa^*)\}
\end{equation}
is a $C^*$-algebra.
\et
The following theorem, taken from \cite{Ln2}, can be thought of as an exact analogue of classical Korovkin's theorem where the convergence is the weak operator convergence.
\bt
Let $\{\Phi_\lambda\}$ be a net of CP-maps on $\mathbb{B}(\mathbb{H}),$ where $\mathbb{H}$ is a separable, complex Hilbert space. If
\begin{equation}\nonumber
\Phi_\lambda(A)\rightarrow A \,\,\textrm{for every A in the set}\,\,\{I,S,S{S^*}\}
\end{equation}
where $S$ is the  unilateral right shift operator on $\mathbb{H}$, then
\begin{equation}\nonumber
\Phi_\lambda(A)\rightarrow A
\end{equation}

for all $A$ in $\mathbb{B}(\mathbb{H})$, where the mode of convergence is the weak operator convergence.
\et
Now we construct a sequence of completely positive maps using the notions in \cite{st1}.
Let $\mathbb{H}$ be a complex separable Hilbert space and let $\{P_n\}$ be a
sequence of orthogonal projections on $\mathbb{H}$ such that
\begin{equation}\label{projection1} \nonumber
\textrm{dim}(P_n(\mathbb{H}))= n <\infty, \,\,\textrm{for each}\,\, n=1,2,3 \ldots
\end{equation}
\begin{equation}\label{projection2} \nonumber
\textrm{and}\,\,\,\mathop {\textrm{lim}}\limits_{n \to \infty } P_n \left( x \right) = x,
 \,\,\textrm{for every x in} \,\,  \mathbb{H}.
\end{equation}

Let $\{U_n\}$ be a sequence of unitary operators on $ \mathbb{H}$ and let
$\mathbb{B}(\mathbb{H})$ be the Banach algebra of all bounded
operators from
$ \mathbb{H}$ to itself. For each $ A \in \mathbb{B}(\mathbb{H})$, consider the
following truncations $A_n=P_nAP_n$ and
$ {\tilde{U_n}}  = P_n U_n $. Then $A_n$ and $ {\tilde{U_n}}$ can be regarded as
$n \times n$ matrices in $M_n \left(\mathbb{C}
\right)$, by restricting their domain to the range of $P_n$. Note that all the matrices
$ {\tilde{U_n}}$ are unitary. For each n, we define the commutative algebra $M_{{\tilde{U_n}}}$
of matrices as follows.
 \begin{equation}\label{Diagonal} \nonumber
 M_{{\tilde{U_n}}}  = \left\{ {A \in M_n \left( \mathbb{C} \right);{\Tilde{U_n}}^*
 A{\Tilde{U_n}} \,\, \textrm{complex diagonal}} \right\}
\end{equation}
Recall that $M_n \left( \mathbb{C} \right)$ is a Hilbert space with the Frobenius norm,
 \begin{equation} \nonumber
\left\| A \right\|_2^2  = \sum\limits_{j,k = 1}^n {\left| {A_{j,k} } \right|^2 }
\end{equation}
induced by the classical Frobenius scalar product,
\begin{equation} \nonumber
\left\langle A,B\right\rangle= \textrm{trace}\,\,(B^*A)
\end{equation}
with trace(.) being the trace of its argument. i.e. the sum
of all its diagonal entries. We observe that
 $M_{{\tilde{U_n}}}$ is a closed convex set in $M_n \left( \mathbb{C}
\right)$ and hence, corresponding to each $A \in M_n \left(
\mathbb{C}
\right)$, there exists a unique
matrix $P_{{\tilde{U_n}}}(A)$ in $M_{{\tilde{U_n}}}$ such that
\begin{equation} \nonumber
\left\| {A - X} \right\|_2^2  \ge \left\| {A - P_{{\tilde{U_n}}}(A) } \right\|_2^2
\,\,\textrm{for every}\,\, X \in M_{{\tilde{U_n}}}.
\end{equation}
We recall the following two lemmas, which reveal some fundamental properties of
the map $P_{{\tilde{U_n}}}$ for each n.
\bl \cite{st1} \label{Fundamental}
With A,B $\in M_n \left(\mathbb{C} \right)$ and $\alpha,\beta$ complex numbers, we have
\begin{equation}\label{Lemma 1}
 P_{{\Tilde{U_n}}}(A)= {\Tilde{U_n}} \sigma \left( {\Tilde{U_n}}^*
 A{\Tilde{U_n}}  \right){{\Tilde{U_n}}}^*
\end{equation}
where $\sigma\left( X\right)$ is the diagonal matrix having $X_{ii}$
as the diagonal elements.
\begin{equation}\label{Lemma 2}
P_{{\Tilde{U_n}}}(\alpha A+\beta B)=\alpha P_{{\Tilde{U_n}}}(A)+
\beta P_{{\Tilde{U_n}}}(B)
\end{equation}
\begin{equation}\label{Lemma 3}
P_{{\Tilde{U_n}}}({A}^*)= {P_{{\Tilde{U_n}}}(A)}^*
\end{equation}
\begin{equation}\label{Lemma 4}
\textrm{Trace}P_{{\Tilde{U_n}}}({A})= \textrm{Trace}(A)
\end{equation}
\begin{equation}\label{Lemma 5}
\left \|P_{{\Tilde{U_n}}}({A})\right \|= 1 \textrm{ (Operator norm)}
\end{equation}
\begin{equation}\label{Lemma 6}
\left \|P_{{\Tilde{U_n}}}({A})\right \|_F= 1\textrm{(Frobenius norm)}
\end{equation}
\begin{equation}\label{Lemma 7}
{\left \|A-P_{{\Tilde{U_n}}}({A})\right \|_F}^2= {\left \|A\right
 \|_F}^2-{\left \|P_{{\Tilde{U_n}}}({A})\right \|_F}^2
\end{equation}
\el
\bl \cite{bto1} \label{positivity}
If A is a Hermitian matrix, then the eigenvalues of $P_{{\Tilde{U_n}}}({A})$
are contained in the closed interval $[\lambda_1(A),\lambda_n(A)]$, where
$\lambda_j(A)$ are the eigenvalues of A arranged in a non decreasing way.
Hence if A is positive definite, then $P_{{\Tilde{U_n}}}({A})$ is positive definite as well.
\el
Now we define a sequence of positive linear maps on $\mathbb{B}(\mathbb{H})$ as follows.
\bd\label{phins}
For each $ A\in \mathbb{B}(\mathbb{H})$,
$\Phi_n:\mathbb{B}(\mathbb{H})\rightarrow \mathbb{B}(\mathbb{H})$ is defined as
\begin{equation} \nonumber
 \Phi_n(A)=P_{{\tilde{U_n}}}(A_n),
\end{equation}
where $P_{{\tilde{U_n}}}(A_n)$ is as in Lemma (\ref{Fundamental}),
for each positive integer $n$.
\ed
We may call $ \Phi_n(A)$, the \textit{pre-conditioners} of A. One of the straightforward but crucial implications of Lemma (\ref{Fundamental}) is the following theorem.
\bt \label{Theorem 1}
The maps $\{\Phi_n\}$ in the Definition (\ref{phins}), is a sequence of completely positive maps on
$\mathbb{B}(\mathbb{H})$ such that
\begin{itemize}
	\item $\left \|\Phi_n\right \|= 1,\,\, \textrm{for each}\,\,n.$
\item  $\Phi_n$ is continuous in the strong topology of operators for each n.
 \item $\Phi_n(I)=I_n$ for each n where I is the identity operator on $\mathbb{H}.$
\end{itemize}
\et
\bp
From Lemma (\ref{positivity}), it follows that $P_{{\tilde{U_n}}}(.) $
 is a positive linear map for each n.
 Since $M_{{\tilde{U_n}}}$ is a commutative Banach algebra,
 $P_{{\tilde{U_n}}}(.)$ is a completely positive map for each n.
  Hence $ {\Phi_n}$ is a completely positive map, since it is the
 composition of CP-maps ($P_{{\tilde{U_n}}}(.)$ and the maps which send A to $P_nAP_n$ and pull back).
 Now continuity in the strong operator topology follows easily from the definition and the remaining
 part of the theorem follows from the following observations.
 \begin{equation}\nonumber
\left\| {\Phi _n } \right\| = \mathop {\sup }\limits_{\left\| A \right\| = 1,A \in
\mathbb{B}(\mathbb{H})}  \left\|{\Phi _n \left( A \right)} \right\|= \mathop {\sup }\limits_{\left\| A \right\| = 1,A \in B\left( H \right)} \left\|
 P_{\Tilde{U_n}}({A_n}) \right\| = 1
\end{equation}
by the identity (\ref{Lemma 5}) in Lemma (\ref{Fundamental}).
The last part of theorem follows easily from the identity (\ref{Lemma 1}) of
Lemma (\ref{Fundamental}).
\ep
In the next section, a few general Korovkin-type theorems  for completely positive maps are proved with respect to various types of clustering of eigenvalues. As an application, the completely positive maps that arises from pre-conditioners discussed above are considered.
\section{Korovkin-type Theorems}
We introduce  different notions of convergence of sequence of
CP-maps in $B(\mathbb{H})$ in a \textit{distributional} sense. We recall the definitions of different notions of convergence for pre-conditioners in \cite{st1}.
To avoid confusion with the classical notion of strong, weak and operator norm convergence,
 we address these by  \textbf{\textit{strong cluster, }}
\textit{\textbf{weak cluster}}, and \textit{\textbf{uniform cluster}}
to mean the strong, weak and uniform convergence respectively used in \cite{st1}.
\bd\label{Stefano}
Let $\{A_n\}$ and $\{B_n\}$ be two sequences of $n\times n$ Hermitian matrices.
We say that $ A_n$ converges to $B_n$ in \textbf{\textit{strong cluster}} if for any $\epsilon>0$,
 there exist integers
$N_{1,\epsilon},N_{2,\epsilon}$ such that all the singular values
  $ \sigma_j( A_n - B_n) $ lie in the interval $\left[0,\epsilon\right)$ except for
at most $N_{1,\epsilon}$ (independent of the size n) eigenvalues for all
$n> N_{2,\epsilon}$. \newline

If the number $ N_{1,\epsilon}$ does not depend on $\epsilon$, we say that  $ A_n $ converges to $ B_n $
in \textit{\textbf{uniform cluster}}. And if $ N_{1,\epsilon}$ depends on $\epsilon, n$ and is of $o(n)$,
we say that $A_n $ converges to $B_n$ in \textit{\textbf{weak cluster}}.

The following powerful lemma is due to Tyrtyshnikov (see Lemma (3.1) in \cite{tyr}).
\bl \label{Tyrtylemma}
Let $\{A_n\}$ and $\{B_n\}$ be two sequences of $n\times n$ Hermitian matrices. If ${\parallel{A_n-B_n}\parallel_F}^2=O(1),$ then we have convergence in strong cluster. If ${\parallel{A_n-B_n}\parallel_F}^2=o(n)$ then the convergence is in weak cluster.
\el
Using the above notions, we introduce the new notions of convergence of CP-maps in $B(\mathbb{H})$.
\bd \label{mnn}
Let $\{\Phi_n\}$ be a sequence of CP-maps in $B(\mathbb{H})$ and $P_n$ be a sequence of projections
on $\mathbb{H}$ with rank n. We say that
$\{\Phi_n(A)\}$ converges to $A$ in the \textbf{\textit{strong distribution sense}},
if the sequence of matrices $\{P_n\Phi_n(A)P_n\}$ converges to $\{P_nAP_n\}$ in strong cluster
as per Definition (\ref{Stefano}). \newline

   Similarly We say that
$\{\Phi_n(A)\}$ converges to $A$ in the \textbf{\textit{weak distribution sense (uniform distribution sense respectively)}},
if the sequence of matrices $\{P_n\Phi_n(A)P_n\}$ converges to $\{P_nAP_n\}$ in weak cluster (uniform cluster respectively)
as per Definition (\ref{Stefano}).
\ed
\br
The above definitions make sense only in the case when A is a self-adjoint operator in
$B(\mathbb{H})$. In the non self-adjoint case, one may have to translate things in to the
language of $\epsilon-$discs instead of intervals. But we are dealing with self-adjoint case
only.
\er
 Consider a sequence of CP-maps $\{\Phi_n\}$
in $B(\mathbb{H})$ with $\left\|\Phi_n\right\| \leq 1.$
By the compactness of
 $CP(\mathbb{B}(\mathbb{H})),$ in the Kadison's B.W topology, $\{\Phi_n\}$ has limit points.
  Let $\Omega$ be the set of all limit points of $\{\Phi_n\}$.
Next we discuss some properties of the limit points  $\Phi $ in $ \Omega$ .
The relation between $\Phi (A)$ and A for $A \in\mathbb{B}(\mathbb{H})$ are considered here.
\bl
Let $ \Phi \in \Omega$ and let $\{\Phi_{n_{ \alpha}}\}$ be a subnet of $\{\Phi_n\}  $ such
that $\Phi_{n_{ \alpha}}$ converges to $\Phi$ in the B.W
topology. Then for each m, the truncations
 $\Phi_{m,n\alpha}(A)=P_m\Phi_{n_{ \alpha}}(A)P_m$ converges to
 $P_m\Phi(A)P_m$ uniformly in norm as ${n_{ \alpha}}\rightarrow \infty$.
 \el
\bp
This follows immediately since $P_m$ is of finite rank and therefore, on range$(P_m),$
weak, strong and operator norm topologies coincides.
\ep
\br
For each $A\in \mathbb{B}(\mathbb{H})$, note that $A_{n_{\alpha}}-\Phi_{n_{\alpha}}(A)$
converges in the strong operator topology to $A- \Phi(A)$.
Hence, $P_mA_{n_{ \alpha}}P_m -P_m\Phi_{n_{ \alpha}}(A)P_m$ converges
to  $P_mAP_m -P_m\Phi(A)P_m$ in the norm topology for each m.
\er
The above observations can be used to deduce the following result.
\bt \label{ncfiniterank}
Let $A\in \mathbb{B}(\mathbb{H})$ be self-adjoint and $\Phi_n(A)$
converges to $A$ in uniform distribution sense as in Definition (\ref{mnn}).
Then $A-\Phi(A)$ is finite rank.
\et
\bp
We use the notation $\#S$ to denote the number of elements in the set S.
By assumption $P_{n_{ \alpha}}\Phi_{n_{ \alpha}}(A)P_{n_{ \alpha}}$ converges
to $A_{n_{ \alpha}}$ in uniform cluster, as in the Definition (\ref{Stefano}). Hence for each
$\epsilon>0$, there exist positive integers $N_{\epsilon}$ and $N$ such that
\begin{equation} \nonumber
\#(\sigma(A_{n_{ \alpha}}- P_{n_{ \alpha}}\Phi_{n_{ \alpha}}(A)P_{n_{ \alpha}})\cap
\mathbb{R} - (-\epsilon,+\epsilon))\leq N,\,\, \textrm{whenever}\,\,{n_{
\alpha}}>N_{\epsilon}.
\end{equation}

Therefore by Cauchy interlacing theorem,
 \begin{equation} \nonumber
\#(\sigma(P_m(A_{n_{ \alpha}}- P_{n_{ \alpha}}\Phi_{n_{ \alpha}}(A)P_{n_{ \alpha}})P_m)\cap
\mathbb{R} - (-\epsilon,+\epsilon))\leq N,\,\, \textrm{if}\,\,{n_{\alpha}}>N_{\epsilon}\,\, \textrm{and}\,\,{n_{\alpha}} \geq m.
\end{equation}
As ${n_{\alpha}}\rightarrow \infty,$ since $P_m(A_{n_{ \alpha}}- P_{U_{n_{ \alpha}}}(A_{n_{ \alpha}}))P_m$
 converges to $P_m(A-\Phi(A))P_m$ in the operator norm topology for every m,
\begin{equation} \nonumber
\#(\sigma(P_m(A-\Phi(A))P_m)\cap \mathbb{R} - (-\epsilon,+\epsilon))\leq N,
\,\, \textrm{for every}\,\, \,m.
\end{equation}
Therefore by Arveson's theorem in \cite{Arv3}, $\mathbb{R} - (-\epsilon,+\epsilon)$ contains no essential points and
hence the essential spectrum $\sigma_e(A- \Phi(A))\subseteq (-\epsilon,+\epsilon) \,\, \textrm{for all}\,\, \epsilon>0.$
This implies that $\sigma_e(A- \Phi(A))=\{0\}.$ Hence $A-\Phi(A)$ is compact and it has at most N
eigenvalues. Hence it is finite rank by spectral theorem.
\ep
\br
Under the assumptions that $A\in \mathbb{B}(\mathbb{H})$ is self-adjoint and
$\Phi_n(A)$ converges to $A$ in uniform distribution sense as in Definition (\ref{mnn}), the
following
results are easy consequences of the above theorem:
\begin{itemize}
\item A is compact if and only if $\Phi(A)$ is compact.
\item A is Fredholm  if and only if $\Phi(A)$ is Fredholm.
\item A is Hilbert Schimidt  if and only if $\Phi(A)$ is Hilbert Schimidt.
\item A is of finite rank if and only if $\Phi(A)$ is of finite rank.
\item A has a gap in the essential spectrum $\sigma_{ess}(A)$ of A if and only if
$\sigma_{ess}(\Phi(A))$ has a gap.
\end{itemize}
\er
\bt \label{nccompact}
Let $A\in \mathbb{B}(\mathbb{H})$ be self adjoint and $\Phi_n(A)$
converges to $A$ in strong distribution sense as in Definition (\ref{mnn}).
 Then $A-\Phi(A)$ is compact.
\et
\bp
The proof is not much different from the proof of Theorem (\ref{ncfiniterank}).
All the arguments are same, except the fact that here the number of eigenvalues of
$(A_{n_{ \alpha}}- P_{n_{ \alpha}}\Phi_{n_{ \alpha}}(A)P_{n_{ \alpha}})$, outside
$(-\epsilon,+\epsilon)$, is not bounded by a constant, but by a
 number $N_{1, \epsilon}$, which depends on $\epsilon$.
Hence we can conclude that $A- \Phi(A)$ is
compact, and can have countably infinite number of eigenvalues.
\ep
\br
Under the assumptions that $A\in \mathbb{B}(\mathbb{H})$ is self-adjoint and
$\Phi_n(A)$ converges to $A$ in strong distribution sense as in Definition (\ref{mnn}), the
following results are easy consequences of the above theorem:
\begin{itemize}
\item A is compact if and only if $\Phi(A)$ is compact.
\item A is Fredholm  if and only if $\Phi(A)$ is Fredholm.
\item A has a gap in the essential spectrum $\sigma_{ess}(A)$ of A if and only if
$\sigma_{ess}(\Phi(A))$ has a gap.
\end{itemize}
\er
Now we recall the notion of the generalized Schwarz map
and a Schwarz type inequality due to Uchiyama \cite{Uchyama}.

\br
Recall Definition (\ref{Schwarz map}), a Schwarz map is clearly positive and contractive. If $\mathbb{A}$ is commutative,
then a positive contractive map is a Schwarz map. A positive linear map $\Psi$	
with $\Psi(1)\leq 1$ was called a Jordan-Schwarz map in \cite{Beck}.
\er
Consider a binary operation $\circ$ in $\mathbb{A}$, satisfying the following
properties for every $\alpha, \beta \in  \mathbb{C}$ and $x, y, z \in \mathbb{A}.$

\begin{enumerate}
	\item \label{1} $ (\alpha x + \beta y)\circ z = \alpha(x \circ z) + \beta(y \circ z).$
	  \item \label{2} $(x \circ y)^* = y^* \circ x^*.$
   \item \label{3} $x^* \circ x \geq 0.$
   \item \label{4} There is a real number M such that $\left\|x \circ y \right\| \leq M\left\|x\right\|\left\|y\right\|.$
 \item \label{5} $(x \circ y) \circ z = x \circ (y \circ z)$
  \item \label{6} $(x \circ y) = (y \circ x)\,$ and $\,x \circ x = x^2\,$ if $\,x = x^*.$
\end{enumerate}
\br
We remark that $\circ$ is bilinear and that the ordinary product satisfies (\ref{5})
and the Jordan product (\ref{6}), conversely if $\circ$ satisfies (\ref{6}), then $\circ$ is the Jordan
product.
\er
\bd \cite{Uchyama}
 A linear map $\Phi$ on $\mathbb{A}$ is called a generalized Schwarz map
with respect to the binary operation $\circ$, if $\Phi$ satisfies
$\Phi(x^*) = \Phi(x)^*$ and $\Phi(x^*) \circ \Phi(x) \leq \Phi(x^* \circ x)$ for every $x \in \mathbb{A}.$
\ed
\br
Note that a generalized Schwarz map $\Phi$ is not necessarily positive. However, under point wise product in function spaces and with usual product of operators or matrices, all Schwarz maps are positive.
\er
\bt \cite{Uchyama} \label{Uchyama}
Let $\Phi$ be a generalized Schwarz map on $C^*$ algebra $\mathbb{A}$, and let

 $\begin{array}{*{20}c}
   {X = } & {\Phi( f^*\circ f)-\Phi(f)^*\circ\Phi(f)}& {\geq0}  \\
   {Y = } & {\Phi( g^*\circ g)-\Phi(g)^*\circ\Phi(g)}& {\geq0}   \\
   {Z = } & {\Phi( f^*\circ g)-\Phi(f)^*\circ\Phi(g)} &&{}  \\
\end{array}\\ $
Then we have for every $f,g \in \mathbb{A},$
\begin{equation} \label{Schwarz}
  \mid \phi(Z)|  \leq \mid\phi(X)|^{1/2}.\mid\phi(Y)|^{1/2},\,\,\textrm{for all state} \,\,
   \phi\,\, on\, \,\mathbb{A}
  \end{equation}
      In particular we have ,
  \begin{equation} \label{Schwarz1}
      \parallel{Z}\parallel \leq \parallel{X} \parallel^{1/2}
      \parallel{Y}\parallel^{1/2}
      \end{equation}
\et
We use these tools to prove the non commutative analogue of the remainder
estimate in the classical Korovkin-type theorems, as proved in \cite{st1}.
 \bl \label{ncremainder}
 Let $\{A_1,A_2,\ldots A_m\}$ be a finite set of operators in  $\mathbb{B}(\mathbb{H})$ and
 $\Phi_n$ be a sequence of  positive linear maps on $\mathbb{B}(\mathbb{H})$
 such that $\left\|\Phi_{n}\right\|\leq1,$ for every n and
 \begin{equation} \nonumber
 \left\|\Phi_{n}(A)-A\right\|=o(\theta_{n})\,\,
 \end{equation}
 for every A in the set $D=\{A_1,A_2,\ldots A_m, \sum_{k=1}^m {A_k  {A_k}^*}\}$,
 where  $\theta_{n}\longrightarrow \,\,0$ as $ n \longrightarrow  \infty.$
 Then  $ \left\|\Phi_{n}(A)-A\right\|=o(\theta_{n})$ for every A in the $C^*$ algebra generated
 by $\{A_1,A_2,\ldots A_m\}.$
  \el
 \bp
We have by linearity,
  \begin{equation} \nonumber
     \Phi_{n}(\sum_{k = 1}^m{A_k {A_k}^*}) = \sum_{k = 1}^m{\Phi_{n}(A_k  {A_k}^*)}
  \end{equation}
Also\\

$\begin{array}{*{20}c}
   {\Phi_{n}(\sum_{k = 1}^m{A_k {A_k}^*})- \sum_{k = 1}^m{A_k  {A_k}^*}}&{ = [\sum_{k =
1}^m\Phi_{n}({A_k} {A_k}^*)-\sum_{k =1}^m{\Phi_{n}(A_k)\Phi_{n}({A_k})^*} ]} & {}  \\\\
   {}&+{ [\sum_{k = 1}^m{\Phi_{n}(A_k)\Phi_{n}(A_k)^*}- {\sum_{k = 1}^m{A_k}  {A_k}^*}]}  & {}  \\
   \end{array}$\\

The norm of left side of the above equation as well as of the last term of the right side
  are of $o(\theta_{n})$. The first term of the right side is
 \begin{equation}\nonumber
\sum_{k = 1}^m[\Phi_{n}(A_k {A_k}^*)-{\Phi_{n}(A_k)\Phi_{n}(A_k)^*}]
\end{equation}
  Hence norm of this term is of $o(\theta_{n})$. But each term inside this sum are
  nonnegative operators by Schwarz inequality for positive linear maps.
 Therefore norm of each term namely $\Phi_{n}(A_k
 {A_k}^*)-{\Phi_{n}(A_k)\Phi_{n}(A_k)^*}$ is of $o(\theta_{n})$. Also since
 each $\Phi_n$ is a Schwarz map, by applying inequality (\ref{Schwarz1}) to maps $\Phi_{n}$ for each n and
 operators $A_k$ and $A_l$, we get
  \begin{equation} \label{ncschwarz}
 \left\|\Phi_{n}({A_k}^*{A_l})-\Phi_{n}({A_k})^*\Phi_{n}({A_l})\right\|=o(\theta_{n}).
  \end{equation}
Also, we can manipulate $\left\|{A_k}^*{A_l}-\Phi_{n}({A_k})^*\Phi_{n}({A_l})\right\|$ as follows.
\begin{eqnarray*}\nonumber
\left\|{A_k}^*{A_l}-\Phi_{n}({A_k})^*\Phi_{n}({A_l})\right\|=\\
\left\|({A_k}-\Phi_{n}({A_k})+\Phi_{n}({A_k}))^*{A_l}-\Phi_{n}({A_k})^*(\Phi_{n}({A_l})-A_l+A_l)\right\|\\
\leq \left\|({A_k}-\Phi_{n}({A_k}))^*{A_l}\right\|+\left\|\Phi_{n}({A_k})^*(\Phi_{n}({A_l})-A_l)\right\|
\end{eqnarray*}
Now each of the terms in the last sum is of $o(\theta_{n})$, since
by the assumption on $A_K,A_l$ and since $\left\|\Phi_{n}\right\|\leq1.$
Therefore we have
  \begin{equation} \label{ncschwarz1}
\left\|{A_k}^*{A_l}-\Phi_{n}({A_k})^*\Phi_{n}({A_l})\right\|= o(\theta_{n})
  \end{equation}
 \begin{equation} \nonumber
\textrm{Now}\,\,\left\|\Phi_{n}({A_k}^*{A_l})-{A_k}^*{A_l}\right\|=\left\|\Phi_{n}({A_k}^*{A_l})-\Phi_{n}({A_k})^*\Phi_{n}({A_l})
+\Phi_{n}({A_k})^*\Phi_{n}({A_l})-{A_k}^*{A_l}\right\|
 \end{equation}
 Applying (\ref{ncschwarz}) and (\ref{ncschwarz1}) in the above identity, we get
  \begin{equation} \nonumber
 \left\|\Phi_{n}({A_k}^*{A_l})-{A_k}^*{A_l} \right\|=o(\theta_{n})
 \end{equation}
 Therefore the proof is completed for every operator of the form ${A_k}^*{A_l}$
 and hence in the algebra generated by the finite set $\{A_1,A_2,\ldots A_m\}$.
 Using the continuity of $\Phi_n's$, we conclude that
 \begin{equation} \nonumber
 \left\|\Phi_{n}(A)-A\right\| =o(\theta_{n})
 \end{equation}
 holds for every operator in the $C^*$ algebra generated by the finite set
$\{A_1,A_2,\ldots A_m\}$.
Hence the proof.
\ep
Before proving the more general Korovkin type theorems, we prove
the following lemma, which is useful for us.
\bl \label{mixed}
 Let $\{A_n\}$ and $\{B_n\}$ be two  sequences of $n\times n$ Hermitian  matrices such that $\{A_n\}$
 converges to $\{B_n\}$ in strong cluster (weak cluster respectively). Assume that $\{B_n\}$ is positive definite and invertible such that
$$
 {B_n} \geq \delta{I_n} > 0,\,\,\textrm{ for all n.}
 $$
Then for a given $\epsilon> 0,$ there will exist positive integers $N_{1,\epsilon},N_{2,\epsilon}$
  such that all eigenvalues of ${B_n}^{-1}{A_n}$ lie in the interval $(1-\epsilon , 1+\epsilon)$
   except possibly for $N_{1,\epsilon}=O(1)$ ($N_{1,\epsilon}=o(n)$ respectively) eigenvalues for every $n>N_{2,\epsilon}$.
\el
 \bp
 First we observe that, since $\{A_n\}$
converges to $\{B_n\}$ in strong cluster (weak cluster respectively),
 by definition, there exists integers $N_{1,\epsilon}, N_{2,\epsilon}$ such that
all eigenvalues of $A_n-B_n$ lie in the interval $(-\epsilon,\epsilon)$
except for at most $N_{1,\epsilon}(N_{1,\epsilon}=o(n)$ respectively)
eigenvalues whenever $n\geq N_{2,\epsilon}$.
 Hence by spectral theorem there exist orthogonal projections $P_n$ and $Q_n$
whose ranges are orthogonal such that
\begin{equation}\nonumber
\textrm{rank}(P_n) +\textrm{rank}(Q_n)=n,\, \textrm{rank}(Q_n) \leq N_{1,\epsilon}, \, \parallel{P_n{(A_n-B_n)}P_n}\parallel<\varepsilon
  \end{equation}
\begin{equation} \label{projsum} \nonumber
 \textrm{and}\,\,\,\,\, A_n-B_n = {P_n (A_n-B_n) P_n + Q_n(A_n-B_n)Q_n }
  \end{equation}
 Hence for $\epsilon\delta>0,$ we have the following decomposition.
  \begin{equation} \label{decomposition}
 {A_n}-{B_n} = R_n+N_n, \,\,\textrm{for all}\,\,n \geq N_{2,\epsilon},
 \end{equation}
where the rank of $R_n$ is bounded above by $N_{1,\epsilon}$ and
  $\left\|N_n\right\|\leq \epsilon\delta.$
     Now let $\beta$ be an eigenvalue of ${B_n}^{-1}{A_n}$ with x being
 the associated eigenvector of norm one. Then we have
   $$
   {B_n}^{-1}{A_n}(x) = \beta{x}.
   $$
   Hence,
   $$
  ({A_n}-{B_n})(x)=(\beta-1){B_n}(x).
   $$
 Which implies that
   $$
   \langle {({A_n}-{B_n})(x),x\rangle =(\beta-1)\langle {B_n}(x),x\rangle}.
   $$
   And
   $$
   {\beta-1} =  \frac{\langle {({A_n}-{B_n})(x),x\rangle}}{\langle {B_n}(x),x\rangle}
   $$
   Now from the decomposition (\ref{decomposition}), we have
      $$
   {\beta-1} =  \frac{\langle ({R_n}+{N_n})(x),x\rangle}{\langle {B_n}(x),x\rangle} =
   \frac{\langle {R_n}(x),x\rangle}{\langle {B_n}(x),x\rangle}+\frac{\langle {N_n}(x),x\rangle}{\langle {B_n}(x),x\rangle}
   $$
 Now since $\left\|N_n\right\|\leq \epsilon\delta,$ and ${B_n} \geq \delta{I_n} > 0,$ the
 second term in the last sum is less than $\epsilon.$
 Also since rank of $R_n$ is bounded above by $N_{1,\epsilon}=O(1)(o(n)$ respectively),
 there are only at most $N_{1,\epsilon}$ linearly independent vectors x for which $R_n(x)\neq 0,$
 by rank-nullity theorem. Hence, except for at most $N_{1,\epsilon}=O(1)(o(n)$ respectively)
 eigenvalues,
  $$
   {\left|\beta-1\right|}\leq \epsilon.
  $$
   This means that all eigenvalues of ${B_n}^{-1}{A_n}$ lie in the interval$(1-\epsilon , 1+\epsilon)$
  except possibly for $N_{1,\epsilon}=O(1)(o(n)$ respectively). This completes the proof.
 \ep
Now we prove our main results, the non commutative versions of Korovkin-type theorems.

\bt \label{nckorovkin}
 Let $\left\{A_1,A_2,\,.\,.\,.A_m\right\}$ be a finite set of self-adjoint
 operators on $\mathbb{H}$ such that $\Phi_n(A)$ converges to $A$
in the strong (or weak respectively) distribution sense, for A in  $\left\{A_1,A_2,\ldots A_m,{A_1}^2,{A_2}^2,\ldots {A_m}^2\right\}$. In addition, if we assume that the difference
 ${P_n(A_k^2)P_n-(P_n(A_k)P_n)^2}$ converges to the 0 matrix in strong cluster (weak cluster respectively),
 for each k, then $\Phi_n(A)$ converges to $A$
in the strong (or weak respectively) distribution sense, for all A in the $C^*$- algebra  $\mathbb{A}$ generated by
 $\left\{A_1,A_2,A_3,\,.\,.\,.\,.A_m\right\}.$
\et
 \bp
  First we consider the following sequence of Hermitian matrices.
 \begin{eqnarray*}
     X_n &=& P_n\Phi_n(A_k^2)P_n- (P_n\Phi_n(A_k)P_n)^2\,\,\geq0  \\
     Y_n &=& P_n\Phi_n(A_l^2)P_n- (P_n\Phi_n(A_l)P_n)^2\,\,\geq0\\
    Z_n &=& P_n\Phi_n(A_k\circ A_l)P_n- (P_n\Phi_n(A_k)P_n)(P_n\Phi_n(A_l)P_n)
\end{eqnarray*}

   (Here $\circ$ denotes the composition of operators).
Since these sequences of matrices are norm bounded, we have
\begin{equation} \label{ncYbound}
\parallel{Y_n}\parallel < \gamma<\infty \,\,\,\textrm{for all n, for some } \gamma>0.
\end{equation}
Also if we write
 \begin{eqnarray*}
X_n &=& P_n\Phi_n(A_k^2)P_n- (P_n\Phi_n(A_k)P_n)^2\\
&=&[P_n\Phi_n(A_k^2)P_n- P_n(A_k^2)P_n]+[P_n(A_k^2)P_n-(P_n(A_k)P_n)^2]\\
&+&[(P_n(A_k)P_n)^2-(P_n\Phi_n(A_k)P_n)^2]
\end{eqnarray*}
the first two terms in the above sum, converges to 0 in strong cluster (weak respectively)
 by assumption. By a simple computation, we get the third term also converges
to the 0 matrix in strong cluster (weak respectively).
 Hence $X_n$ converges to the 0 matrix in strong cluster (weak respectively).

Now for each fixed x with $\left\|x\right\|=1,$ if we consider the state $\phi_x$
on $B(\mathbb{H})$ defined as
\begin{equation}\nonumber
\phi_x(A)=\left\langle A(x),x\right\rangle
\end{equation}
then by the inequality (\ref{Schwarz}), we get
\begin{equation} \label{minmax}
\left|\left\langle Z_n(x),x\right\rangle\right| \leq {\left|{\left\langle X_n(x),x\right\rangle}\right|}^{1/2}.{\left|{\left\langle Y_n(x),x\right\rangle}\right|}^{1/2}
\end{equation}

 Now let $\delta>0,$ be given and $\epsilon={{\delta}^2}/\gamma$,
as in the proof of Lemma (\ref{mixed}), there exists integers  $N_{1,\epsilon}=O(1)$
$(N_{1,\epsilon}=o(n)$ respectively), $N_{2,\epsilon}$
such that we have the following decomposition
\begin{equation}\nonumber
X_n=N_n+R_n \textrm{ for all }n>N_{2,\epsilon},
\end{equation}
with $\left\|N_n\right\|<\epsilon$ and rank of $R_n$ is less than $N_{1,\epsilon}=O(1)$
$(N_{1,\epsilon}=o(n)$ respectively).
Applying this in the inequality (\ref{minmax}), we get
\begin{equation} \nonumber
\left|\left\langle Z_n(x),x\right\rangle\right| \leq \sqrt{\gamma}.[{\left|{\left\langle N_n(x),x\right\rangle}\right|}^{1/2}+{\left|{\left\langle R_n(x),x\right\rangle}\right|}^{1/2}]
\,(\textrm{for all } n>N_{2,\epsilon})
\end{equation}
 Since the rank of $R_n$ is bounded above by $N_{1,\epsilon}=O(1)(o(n)$ respectively),
 there are only at most $N_{1,\epsilon}$ linearly independent vectors x for which $R_n(x)\neq 0,$
 by rank-nullity theorem. Hence,
 $\left|\left\langle Z_n(x),x\right\rangle\right| \leq \delta,$
  except for at most $N_{1,\epsilon}=O(1)(o(n)$ respectively) linearly independent vectors x.
Therefore all eigenvalues of $Z_n$, except for possibly $N_{1,\epsilon}=O(1)(o(n)$ respectively),
lie in the interval $(-\delta, \delta)$, whenever $n>N_{2,\epsilon}.$
 Since $\delta>0,$ was arbitrary, $Z_n$ converges to the 0 matrix in strong cluster (weak respectively).

Now consider
 \begin{eqnarray*}
{P_n\Phi_n(A_k\circ A_l)P_n-P_nA_k\circ A_lP_n}&=& [P_n\Phi_n(A_k\circ A_l)P_n-(P_n\Phi_n(A_k)P_n)(P_n\Phi_n(A_l)P_n)]\\
    &+&[(P_n\Phi_n(A_k)P_n)(P_n\Phi_n(A_l)P_n)-(P_nA_kP_n)(P_nA_lP_n)]\\
    &+&[(P_nA_kP_n)(P_nA_lP_n)-P_n(A_k\circ A_l)P_n]
 \end{eqnarray*}
 By similar arguments above, we see that each term in the right hand side of the above equation converges to the 0 matrix in strong cluster (weak respectively). Hence the theorem is proved for the operators of the form $A_k\circ A_l$. Hence it is true for any operator in the algebra generated by  $\left\{A_1,A_2,A_3,\,.\,.\,.\,.A_m\right\}.$ \\

 Now for $A\in \mathbb{A},\,\epsilon > 0$, let T be the operator in the algebra generated by  $\left\{A_1,A_2,A_3,\,.\,.\,.\,.A_m\right\},$ such that
 \begin{equation} \nonumber
  \parallel{A-T}\parallel < \epsilon/3,\, \textrm{and}
\parallel{\Phi_n(A)-\Phi_n(T)}\parallel<\epsilon/3.
 \end{equation}
 Consider the following equation:\\
 \begin{eqnarray*}
   P_n\Phi_n(A)P_n-P_nAP_n&=&[P_n\Phi_n(A)P_n-P_n\Phi_n(T)P_n]+[P_n\Phi_n(T)P_n-P_nTP_n]\\
   &+& [P_nTP_n-P_nAP_n]
 \end{eqnarray*}

Thus the norm of the sum of the first and third terms is less than $2\epsilon/3.$
The middle term $P_n\Phi_n(T)P_n-P_nTP_n$ can be split into a term with norm less than
$\epsilon/3$ and a term with constant rank independent of the order $n$
(or of o(n) respectively) since T is in the
algebra generated by $\left\{A_1,A_2,A_3,\,.\,.\,.\,.A_m\right\}$.
Thus the sequence of matrices $P_n\Phi_n(A)P_n$ converges to $P_nAP_n$ in
strong cluster (or in weak cluster respectively).
Hence the proof is completed.
\ep
\section{Applications to matrix algebras}
In this section, we will be dealing with the completely positive maps $\Phi_n$,
that we introduced in Definition (\ref{phins}), and its modifications.
We treat them as examples for the results in last section.

Recall the Definition (\ref{Stefano}), given a sequence of algebras $M_{{\tilde{U_n}}}$, with associated operator $P_{{\Tilde{U_n}}}(.)$,
we say that $P_{{\Tilde{U_n}}}(A_n)$ converges to $A_n$ in strong cluster, if for any $\epsilon>0$,
 there exist integers $N_{1,\epsilon},N_{2,\epsilon}$ such that all the singular values
$\sigma_j( A_n - P_{\tilde{U_n}}(A_n) )$
lie in the interval $\left[0,\epsilon\right)$ except for at most $N_{1,\epsilon}$
(independent of the size n) eigenvalues for all $n> N_{2,\epsilon}$. If the number $N_{1,\epsilon}$ does not depend on $\epsilon$, we say that $P_{{\Tilde{U_n}}}(A_n)$ converges to $A_n$ in uniform cluster.
 And if $N_{1,\epsilon}$ depends on $\epsilon, n$ and of $o(n)$, we say that
$P_{{\Tilde{U_n}}}(A_n)$ converges to $A_n$ in weak cluster.
\ed
 Consider the sequence $\{\Phi_n\}$ of Definition (\ref{phins}). By the compactness of
 $CP(\mathbb{B}(\mathbb{H})),$ in the Kadison's B.W topology, $\{\Phi_n\}$ has limit points.
We note some of the properties of the limit points $\Phi$ of $\{\Phi_n\}$, as immediate
consequences of Theorem (\ref{ncfiniterank}) and Theorem (\ref{nccompact}).
 The following  special cases are of interest spectral theory point of view.
\begin{itemize}
\item A is a Hilbert Schmidt operator on $\mathbb{H}$.
\item $A=A(f)$, is the Toeplitz operator where the symbol function $f \in C\left[-\pi,\pi\right]$ and
$\mathbb{H}=L^2\left[-\pi,\pi\right]$.
\item A is a Fredholm or compact operator.
\end{itemize}
\bt \label{finiterank}
Let $A\in \mathbb{B}(\mathbb{H})$ be self adjoint and $P_{{\Tilde{U_n}}}(A_n)$
converges to $A_n$ in uniform cluster  as in Definition (\ref{Stefano}). Then $A-\Phi(A)$ is finite rank.
\et
\bp
Follows easily from Theorem (\ref{ncfiniterank}), by considering
$\Phi_n(A)=P_{{\Tilde{U_n}}}(A_n)$.
\ep
\bt \label{compact}
Let $A\in \mathbb{B}(\mathbb{H})$ be self-adjoint and $P_{{\Tilde{U_n}}}(A_n)$
converges to $A_n$ in strong cluster as in Definition \ref{Stefano}. Then $A-\Phi(A)$ is compact.
\et
\bp
Follows easily from Theorem (\ref{ncfiniterank}), by considering
$\Phi_n(A)=P_{{\Tilde{U_n}}}(A_n)$.
\ep
\br
The analysis of convergence in `weak cluster', in the sense of Def (\ref{Stefano}) is taken up
later in this article.
\er
\subsection{Modified pre-conditioners}
It is interesting to observe that the notion of pre-conditioners can be modified by replacing '\textit{diagonal transformation'} by '\textit{pinching functions'} as follows.\\

 Let $M_{{U_n}}$  = $\left\{ {A \in M_n \left( \mathbb{C} \right);{U_n}^* A{U_n} \ \
 \textrm{is block diagonal}} \right\}$, where the
 block diagonal is  obtained for each $A$ in $M_n(C)$ by applying pinching function
 to $A$ for each n (see \cite{Bha1} for definition). To be more precise, let
 $P_{n_k}$
 be  pairwise orthogonal, $m_n$ orthogonal projections in $M_n(C),$ such that
 $\sum_{k=1}^{m_n} P_{n_k}$ =${I_n}$, the identity matrix. The
 modified
 pre-conditioner on $M_n(C)$ takes values
\begin{equation}\label{Block}
 \Psi_n(A)=\sum_{k=1}^{m_n} P_{n_k}A P_{n_k}\, \textrm{for every}\, A \in M_n(C).
\end{equation}

   From Stinespring's theorem, it is clear that the maps $\Psi_n'$s are CP-maps. Now
 if we define $P_{{U_n}}(A)$ in a similar way
 with $M_{\tilde{U_n}}$ replaced by $M_{{U_n}}$,
  we can formulate an analogue of Lemma (\ref{Fundamental}).
   \bl \label{Fundamentalb}
With A,B $\in M_n \left(\mathbb{C} \right)$, we have
\begin{equation}\label{Lemma b1} \nonumber
P_{{U_n}}(A)= {U_n}\Psi_n \left( {U_n}^* A{U_n}  \right){{U_n}}^* \,\,\textrm{where}\,\,
\Psi_n \,\,\textrm{is as in}\,\,(\ref{Block}).
\end{equation}
\begin{equation}\label{Lemma b2}  \nonumber
 P_{{U_n}}(\alpha A+\beta B)=\alpha P_{{U_n}}(A)+ \beta P_{{U_n}}(B)
\end{equation}
\begin{equation}\label{Lemma b3}  \nonumber
P_{{U_n}}({A}^*)= {P_{{U_n}}(A)}^*
\end{equation}
\begin{equation}\label{Lemma b4}  \nonumber
Trace P_{{U_n}}({A})= Trace (A)
\end{equation}
\begin{equation}\label{Lemma b5}  \nonumber
 \left \|P_{{U_n}}({A})\right \|= 1 \textrm{ (Operator norm)}
\end{equation}
\begin{equation}\label{Lemma b6}  \nonumber
\left \|P_{{U_n}}({A})\right \|_F= 1\textrm{(Frobenius norm)}
\end{equation}
\begin{equation}\label{Lemma b7}  \nonumber
{\left \|A-P_{{U_n}}({A})\right \|_F}^2= {\left \|A\right
\|_F}^2-{\left \|P_{{U_n}}({A})\right \|_F}^2
\end{equation}
\el
We shall list down some of the properties of the maps $\{\Psi_n\}$ as we did in
Theorem (\ref{Theorem 1}).
 \bt
The maps $\{\Psi_n\}$ is a sequence of completely positive maps on
$\mathbb{B}(\mathbb{H})$ such that
\begin{itemize}

        \item \label{Theorem b1}  \nonumber
                 $\left \|\Psi_n\right \|= 1,$ for each n.
        \item $ \Psi_n$ is continuous in the strong topology of operators.
        \item $\Psi_n(I)=I_n$ for each n, where I is the identity operator on
        $\mathbb{H}.$
\end{itemize}
\et
 \br
The above mentioned modified version of pre-conditioners is better than
the previous one in the sense that the modified version is
closer to the operator in the Frobenius norm and is simpler
enough also.
 \er
 We may construct examples for the modified pre-conditioners as follows.
 \be \label{vander}
 Let $ {\tilde{U_n}}$ be unitaries in $M_n(C)$
as in Definition (\ref{Stefano}). For each positive integer $n$,
let ${U_n}$ be unitaries in $\mathbb{B}(\mathbb{H})$ defined as
$\tilde{U_n}\bigoplus(I-{P_n})$. Observe that there are many
interesting, concrete examples of unitaries ${\tilde{U_n}}$ in $\cite{st1}$.
For the sake of completeness, we quote them below.

Let $v = \left\{ {v_n } \right\}_{n \in N} \,\, \textrm{with}\,v_n  = \left(
{v_{nj}
} \right)_{j \le n - 1}$ be a sequence of trigonometric
functions on an interval I. Let $S = \left\{ {S_n } \right\}_{n \in N}$
be a sequence of grids of n points on I,
namely, $S_n  = \left\{ {x_i^n ,i = 0,1,.\,.\,.\,n - 1} \right\}$.
Let us suppose that the generalized Vandermonde matrix
\begin{equation} \nonumber
{ V_n  = }\left( {{ v}_{{ nj}} { (}x_i^n { )}} \right)_{{ i;j = 0}}^{{ n - 1}}
\end{equation}
is a unitary matrix. Then, algebra of the form $M_{{U_n}}$ is a trigonometric
algebra if $U_n = {V_n}^{*}$ with $V_n$ a generalized trigonometric Vandermonde matrix.

We get examples of trigonometric algebras with the following choice of the
sequence of matrices ${U_n}$ and grid $S_n$.

\begin{eqnarray*}
U_n &=& F_n=\left( {\frac{{ 1}}{{\sqrt { n} }}e^{ijx_i^n } } \right)\,\,,\,\,\,\,i,j
=0,1,.\,.\,.\,n - 1,\\
S_n&=& \left\{ {x_i^n  = \frac{{2i\pi }}{n},i = 0,1,.\,.\,.\,n - 1} \right\} \subset
I= \left[ { - \pi ,\pi } \right]\\
U_n &=& G_n =\left( {\sqrt {\frac{2}{n+1}} }sin(j+1)x_i^n \right)\,\,,\,\,\,\,i,j =
0,1,.\,.\,.\,n - 1,\\
S_n&=& \left\{ {x_i^n  = \frac{{(i+1)\pi }}{n+1},i = 0,1,.\,.\,.\,n - 1} \right\}
\subset I = \left[ {0 ,\pi } \right]\\
U_n&=& H_n =\left( \frac{{1}}{{\sqrt {n} }}\left[ {{sin(jx_i^n) + cos(jx_i^n )}}
\right] \right)\,\,,\,\,\,\,i,j = 0,1,.\,.\,.\,n - 1,\\
S_n&=&\left\{ {x_i^n  = \frac{{2i\pi }}{n},i = 0,1,.\,.\,.\,n - 1} \right\} \subset
I= \left[ { - \pi ,\pi } \right]
\end{eqnarray*}
\ee
\subsection{Korovkin type theory for Toeplitz operators}
Now we consider the case where $A=A(f)$, is the Toeplitz operator where
the symbol function $f \in C\left[-\pi,\pi\right]$ and $\mathbb{H}=L^2\left[-\pi,\pi\right]$.
We generalize some of the results in \cite{st1} and get stronger versions.
First we recall the Korovkin type results in \cite{st1}.
 We use the notation
$A_n(f)$ for the finite Toeplitz matrix with symbol f.
\bt \label{Stefano1}\cite{st1} Let f be a continuous periodic real-valued function.
Then $P_{U_n}(A_n(f))$
converges to $A_n(f)$ in strong cluster, if $P_{U_n}(A_n(p))$ converges to $A_n(p)$ in strong cluster
for all the trigonometric polynomials p.
\et
\bt \label{Stefano2}\cite{st1} Let f be a continuous periodic real-valued function.
Then $P_{U_n}(A_n(f))$ converges to $A_n(f)$ in weak cluster if $P_{U_n}(A_n(p))$ converges to $A_n(p)$ for all the trigonometric polynomials p.
\et
Before proving the general results, we prove the following lemma, the remainder estimate
version of classical Korovkin's theorem as proved in $\cite{st1}$,
which is used to get more general versions of Theorems (\ref{Stefano1}) and
(\ref{Stefano2}). This is the commutative version of Lemma (\ref{ncremainder}).
 \bl \label{remainder}
 Let $\{g_1,g_2,\ldots g_m\}$ be a finite set of continuous periodic functions and
 $\Phi_n$ be a sequence of  positive linear maps on $C [0,2{\pi}]$
 such that $\left\|\Phi_{n}\right\|\leq1,$ for every n, and
 \begin{equation} \nonumber
 \Phi_{n}(g) =g+o(\theta_{n})\,\,\textrm{for every g in the set} \,\,D=\{g_1,g_2,\ldots g_m
 ,\sum_{k=1}^m {g_k  {g_k}^*}\},
 \end{equation}
 where  $\theta_{n}\longrightarrow \,\,0$ as
 $ n \longrightarrow  \infty.$
 Then  $\Phi_{n}(g) =g+o(\theta_{n})\,\,$ for every g in the $C^*$ algebra generated
 by $\{g_1,g_2,\ldots g_m\}.$
  \el
 \bp
 The proof is obtained by replacing functions in place of operators, in the proof
of Lemma (\ref{ncremainder}).
Using linearity of $\Phi_n$'s, we write,

$
\\ \begin{array}{*{20}c}
   {\Phi_{n}(\sum_{k = 1}^m{g_k {g_k}^*})- \sum_{k = 1}^m{g_k  {g_k}^*}}&{ = \sum_{k =
1}^m\Phi_{n}({g_k} {g_k}^*)-\sum_{k =1}^m{\Phi_{n}(g_k)\Phi_{n}({g_k})^*} } & {}  \\\\
   {}&+{ (\sum_{k = 1}^m{\Phi_{n}(g_k)\Phi_{n}(g_k)^*}- {\sum_{k = 1}^m{g_k}  {g_k}^*})}  & {}  \\
   \end{array}
$\\

The left side of the above equation as well as the last term of the right side
  are of $o(\theta_{n})$. Hence the first term of the right side
 \begin{equation}\nonumber
\sum_{k = 1}^n[\Phi_{n}(g_k {g_k}^*)-{\Phi_{n}(g_k)\Phi_{n}(g_k)^*}]
\end{equation}
 is of $o(\theta_{n})$. But each term inside this sum is
  nonnegative by Schwarz inequality for positive linear maps.
 Therefore each of its terms, namely $\Phi_{n}(g_k
 {g_k}^*)-{\Phi_{n}(g_k)\Phi_{n}(g_k)^*}$ is of $o(\theta_{n})$. Also since every
 positive contractive map in a commutative $C^*$ algebra is a Schwarz map, each $\Phi_n$
 is a Schwarz map.
 Therefore by applying inequality (\ref{Schwarz1}) to the maps $\Phi_{n}$ for each n and
 functions $g_k,g_l$, we get
  \begin{equation} \nonumber
\Phi_{n}({g_k}^*{g_l})-\Phi_{n}({g_k})^*\Phi_{n}({g_l})=o(\theta_{n})
  \end{equation}
Also, we observe the following
 \begin{equation} \nonumber
\Phi_{n}({g_k})^*\Phi_{n}({g_l})-{g_k}^*{g_l}=({g_k}^*+o(\theta_{n}))({g_l}+o(\theta_{n}))-{g_k}^*{g_l} =o(\theta_{n})
 \end{equation}
Using the above two identities, we deduce that
 \begin{equation} \nonumber
\Phi_{n}({g_k}^*{g_l})-{g_k}^*{g_l}=[\Phi_{n}({g_k}^*{g_l})-\Phi_{n}({g_k})^*\Phi_{n}({g_l})]+[\Phi_{n}({g_k})^*\Phi_{n}({g_l})-{g_k}^*{g_l}]=o(\theta_{n})
 \end{equation}
 Therefore the proof is completed for every function of the form ${g_k}^*{g_l}$
 and hence in the algebra generated by $\{g_1,g_2,\ldots g_m\}$.
 Using the continuity of $\Phi_n's$, we conclude that
 \begin{equation} \nonumber
 \Phi_{n}(g)-g =o(\theta_{n})
 \end{equation}
 holds for every function in the $C^*$ algebra generated by the finite set
 $\{g_1,g_2,\ldots g_m\}$.
Hence the proof.
\ep

Now we prove some general versions of Theorems (\ref{Stefano1}) and
(\ref{Stefano2}). The technique of the proof is the same as in Theorem (\ref{nckorovkin}).
Still we provide all the details.
\bt \label{strongkorovkin}
 Let $\left\{g_1,g_2,\,.\,.\,.g_m\right\}$ be a finite set of real valued
 continuous $2\pi$ periodic functions such that $P_{U_n}(A_n(f))$ converges to $A_n(f)$
in strong cluster, for f in  $\left\{g_1,g_2,\ldots g_m,{g_1}^2,{g_2}^2,\ldots {g_m}^2\right\}$.
Then $P_{U_n}(A_n(f))$ converges to $A_n(f)$ in strong cluster for
 all f in the $C^*$- algebra  $\mathbb{A}$ generated by
 $\left\{g_1,g_2,g_3,\,.\,.\,.\,.g_m\right\}.$
\et
 \bp
For any $k,l=1,2,3 \ldots m,$ setting
 \begin{eqnarray*}
   X_n&=& P_{U_n}(A_n(g_k^2))- P_{U_n}(A_n(g_k))^2 \,\,\geq0  \\
     Y_n &=& P_{U_n}(A_n(g_l^2))- P_{U_n}(A_n(g_l))^2 \,\,\geq0\\
    Z_n &=&  P_{U_n}(A_n(g_k^*\circ g_l))- P_{U_n}(A_n(g_k))^*\circ P_{U_n}(A_n(g_l))
\end{eqnarray*}
   (Here $\circ$ denotes the usual point wise product in the case of scalar valued functions and
   matrix product in the case of matrices),
we observe that $ X_n,Y_n,$ and $Z_n$ are all Hermitian matrices of order $n$.
It is clear that all the above sequences of matrices are norm bounded. Then for all
$n$
\begin{equation} \label{Ybound}
\parallel{Y_n}\parallel < \gamma<\infty
\end{equation}
Also if we write
\begin{equation} \nonumber
X_n = \Phi_n(g_k^2)- \Phi_n(g_k)^2=[\Phi_n(g_k^2)- A_n(g_k^2)]+ [A_n(g_k^2)- A_n(g_k)^2]+[A_n(g_k)^2-\Phi_n(g_k)^2],
\end{equation}
the first term on the right hand side of the above equality converges to 0 in strong cluster
 by assumption. The second term is
 \begin{eqnarray*}
   {A_n(g_k^2)- A_n(g_k)^2}&=&{P_nA(g_k^2)P_n- (P_nA(g_k)P_n)^2 } \\
   &=& {P_nA(g_k^2)P_n-P_nA(g_k)^2P_n+P_nA(g_k)^2P_n-(P_nA(g_k)P_n)^2}\\
   &=&{P_n[A(g_k^2)- A(g_k)^2]P_n+P_nA(g_k)^2P_n-(P_nA(g_k)P_n)^2}\\
    &=&{P_n[A(g_k^2)- A(g_k)^2]P_n+P_nH(g_k)^2P_n+Q_nH(g_k)^2Q_n}
\end{eqnarray*}
The last equality is due to Widom \cite{Widom76}, where $Q_n$'s are
projections and $H(g_k)$ is the Hankel operator, which is compact, since the symbols are continuous. Also
 $A(g_k^2)- A(g_k)^2$ is a compact operator
 (for eg. see \cite{Botband}). Hence
 $A_n(g_k^2)- A_n(g_k)^2$ can be written as the sum of sequences of matrices
that are truncations of compact operators. But the compliment
of any neighborhood of 0 contains only finitely many eigenvalues of a compact
operator, as its spectral values. Also the truncations of a compact operator
on a separable Hilbert space converges to the operator in norm.  Therefore
we conclude that $A_n(g_k^2)- A_n(g_k)^2$ converges to the 0 matrix in strong cluster.
By a simple computation, we get that the third term also converges to the 0 matrix in strong cluster.
Hence $X_n$ converges to the 0 matrix in strong cluster.

 By the similar arguments in the proof of Theorem (\ref{nckorovkin}), we
 conclude that $Z_n$ converges to the 0 matrix in strong cluster.

Now consider
 \begin{eqnarray*}
P_{U_n}(A_n(g_k \circ g_l))-A_n(g_k \circ g_l)&=& [P_{U_n}(A_n(g_k \circ g_l))-P_{U_n}(A_n(g_k ))P_{U_n}(A_n(g_l))]\\
   &+&[P_{U_n}(A_n(g_k ))P_{U_n}(A_n(g_l))-A_n(g_k )A_n(g_l)]\\
   &+&[A_n(g_k )A_n(g_l)-A_n(g_k \circ g_l)]
\end{eqnarray*}
By similar arguments above, we see that each term in the right hand side of the above equation converges to the 0 matrix in strong cluster. Hence the theorem is proved for the functions of the form $g_kg_l$. Hence it is true for any function in the algebra generated by $\left\{g_1,g_2,g_3,\,.\,.\,.\,.g_m\right\}.$ \\

 Now for $f \in \mathbb{A},\,\epsilon > 0$, g be the function in the algebra generated by $\left\{g_1,g_2,g_3,\,.\,.\,.\,.g_m\right\}$ such that
 \begin{equation} \nonumber
  \parallel{A_n(f)-A_n(g)}\parallel < \epsilon/3,\, \textrm{and}
\parallel{P_{U_n}(A_n(g))-P_{U_n}(A_n(f))}\parallel<\epsilon/3.
 \end{equation}
 Consider the following equation:
 \begin{eqnarray*}
   A_n(f)-P_{U_n}(A_n(f))&=& [A_n(f)-A_n(g)]+[A_n(g)-P_{U_n}(A_n(g))] \\
   &+& [P_{U_n}(A_n(g))-P_{U_n}(A_n(f))]
\end{eqnarray*}

Thus the norm of the sum of the first and third terms is less than $2\epsilon/3.$
The middle term $A_n(g)-P_{U_n}(A_n(g))$ can be split into a term with norm less than
$\epsilon/3$ and a term with constant rank {independent of the order $n$} since g is in the
algebra generated by $\left\{g_1,g_2,g_3,\,.\,.\,.\,.g_m\right\}$.
Hence the proof is completed.
\ep
\bc
 If $P_{U_n}(A_n(f))$ converges to $A_n(f)$ in strong cluster for all f in $\{1,x,x^2\}$,
 then $P_{U_n}(A_n(f))$ converges to $A_n(f)$ in strong cluster for
 all f in $C[0,2\pi].$
\ec
\bc
 Under the assumption of Theorem (\ref{strongkorovkin}), if
 $f \in \mathbb{A}$ is strictly positive, then for any $\epsilon > 0$, for n
large enough, the matrix\\
 $ P_{U_n}(A_n(f))^{-1}(A_n(f))$ has
eigenvalues in $(1 - \epsilon,1 + \epsilon)$ except for
$N_{\epsilon} = O(1)$ outliers, at most.
\ec
\bp
Since $f \in \mathbb{A}$ is strictly positive, $(A_n(f))$ is positive definite.
This implies that $P_{U_n}(A_n(f))$ is positive definite. Hence the proof is completed
by Lemma (\ref{mixed}).
\ep

Now we prove the exact analogue of Theorem (\ref{strongkorovkin}) in the case of
convergence in weak cluster.
The proof is more or less is the same but for some obvious modifications. However all the details are provided.
\bt \label{weakkorovkin}
 Let $\left\{g_1,g_2,\,.\,.\,.g_m\right\}$ be a finite set of real valued
 continuous $2\pi$ periodic functions such that $P_{U_n}(A_n(f))$ converges to $A_n(f)$ in weak cluster
 for f in  $\left\{g_1,g_2,\ldots g_m,{g_1}^2,{g_2}^2,\ldots {g_m}^2\right\}$.
 Then $P_{U_n}(A_n(f))$ converges to $A_n(f)$ in weak cluster for
 all f in the $C^*$- algebra  $\mathbb{A}$ generated by
 $\left\{g_1,g_2,g_3,\,.\,.\,.\,.g_m\right\}.$
\et
\bp
The proof is the same as Theorem (\ref{strongkorovkin}), except that the
splitting of terms must be as the sum of one with small norm and the other of rank o(n). We give the details below.
Applying (\ref{Schwarz1}) with
 $\Phi_n=P_{U_n}(A_n(.))$ and $X_n,Y_n,Z_n$ as in the proof of Theorem (\ref{strongkorovkin}),
if we write
\begin{eqnarray*}
X_n = \Phi_n(g_k^2)- \Phi_n(g_k)^2&=&[\Phi_n(g_k^2)- A_n(g_k^2)]+ [A_n(g_k^2)- A_n(g_k)^2]\\
&+&[A_n(g_k)^2-\Phi_n(g_k)^2]
\end{eqnarray*}
the first term on the right hand side of the above equality converges to 0 in weak cluster by assumption.
The second term, $A_n(g_k^2)- A_n(g_k)^2$ converges to the 0 matrix in strong cluster
by the same argument in the proof of Theorem (\ref{strongkorovkin}), and hence it converges in weak cluster.
By a simple computation, we get that the third term also converges to the 0 matrix in weak cluster.
Hence $X_n$ converges to the 0 matrix in weak cluster.

By a similar arguments in the proof of Theorem (\ref{nckorovkin}),
we conclude that $Z_n$ converges to the 0 matrix in weak cluster.

Now consider
\begin{eqnarray*}
   P_{U_n}(A_n(g_k \circ g_l))-A_n(g_k \circ g_l)&=& [P_{U_n}(A_n(g_k \circ g_l))-P_{U_n}(A_n(g_k ))P_{U_n}(A_n(g_l))]\\
  &+&[P_{U_n}(A_n(g_k ))P_{U_n}(A_n(g_l))-A_n(g_k )A_n(g_l)]\\
  &+& [A_n(g_k )A_n(g_l)-A_n(g_k \circ g_l)]
\end{eqnarray*}
By similar arguments above, we see that each term in the right hand side of the above equation converges to the 0 matrix in weak cluster. Hence the theorem is proved for the functions of the form $g_kg_l$. Hence it is true for any function in the algebra generated by $\left\{g_1,g_2,g_3,\,.\,.\,.\,.g_m\right\}.$ \\

 Now for $f \in \mathbb{A},\,\epsilon > 0$, g be the function in the algebra generated by $\left\{g_1,g_2,g_3,\,.\,.\,.\,.g_m\right\}$ such that
 \begin{equation} \nonumber
  \parallel{A_n(f)-A_n(g)}\parallel < \epsilon/3,\, \textrm{and}
\parallel{P_{U_n}(A_n(g))-P_{U_n}(A_n(f))}\parallel<\epsilon/3.
 \end{equation}
 Consider the following equation:
\begin{eqnarray*}
   A_n(f)-P_{U_n}(A_n(f))&=&[A_n(f)-A_n(g)]+[A_n(g)-P_{U_n}(A_n(g))] \\
&+&[P_{U_n}(A_n(g))-P_{U_n}(A_n(f))]
\end{eqnarray*}

Thus the norm of the sum of the first and third terms is less than $2\epsilon/3.$
The middle term $A_n(g)-P_{U_n}(A_n(g))$ can be split into a term with norm less than
$\epsilon/3$ and a term with rank $o(n)$, since g is in the
algebra generated by $\left\{g_1,g_2,g_3,\,.\,.\,.\,.g_m\right\}$.
Hence the proof is completed.
\ep
  \bc With the hypotheses Theorem (\ref{weakkorovkin}), if  $f \in \mathbb{A}$ is positive,
then for any $\epsilon > 0,$ for n large enough, the matrix
$P_{U_n}(A_n(f))^{-1}(A_n(f))$ has eigenvalues in $(1 -\epsilon,1 + \epsilon)$
 except $N_{\epsilon} = o(n)$ outliers, at most.
 \ec
 \bp
 Proof follows easily from Lemma (\ref{mixed}).
 \ep
\br
   It is to be noted that Theorem (\ref{strongkorovkin}), (\ref{weakkorovkin}) and
   the corollaries are much stronger than the corresponding theorems in \cite{st1},
    where it has been assumed that the convergence takes place on the algebra
    generated by the \textit{test set}. But here it is assumed that the convergence takes place only on the \textit{test set} as in the classical Korovkin-type theorems. However it is not clear whether the  assumption of convergence on ${g_k}^2$ for each $k$ can be replaced by convergence on $\sum_{k = 1}^n{g_k}^2$ as in the usual case.
 \er
 \section{The LPO sequences }
  It can also be observed that similar stronger versions of Theorems (5.3) and (5.4) of \cite{st1} are
 valid as above. First we recall some of the
 preliminaries from $\cite{st1}$ needed subsequently.\\

 The behavior of eigenvalues of $P_{U_n}(A_n(f))$ has been studied in
 \cite{st1} when ${U_n}$ is the sequence of generalized Vandermonde matrices
 (Example \ref{vander}). Recall that the $j$th row of ${U_n}$ is a vector of trigonometric functions calculated
 on the grid point ${x_j}^{(n)}.$ From Lemma (\ref{Fundamental}), it follows that the $j$th eigenvalue
  $\lambda_j$ of $P_{U_n}(A_n(f))$ is $\sigma(U_nA_n(f){U_n}^*)_{j,j}.$
  Thus $\lambda_j$ is the value of the trigonometric function that takes on
  the $j$th grid point $x={x_j}^{(n)}.$ Now we consider the function $[L_n[U_n](f)](x)$
 obtained by replacing ${x_j}^{(n)}$ by $x$ in $[0,2\pi]$ in the expression of $\lambda_j$.
To make it precise, let $v(x)$ denote the trigonometric function whose values at grid points
$\{{x_j}^{(n)}\},$ form the $j$th generic row of ${U_n}^*$.
 We define the linear operator $L_n[U_n]$ on $C[0,2\pi]$ as follows;
 \begin{equation}
 L_n[U_n](f)= v(x){A_n(f)}v^*(x).
  \end{equation}
  $L_n[U_n](f)$ is the continuous expressions of the diagonal elements of
   $U_nA_n(f){U_n}^*.$
 And it is clear that $L_n[U_n]$ is a sequence of completely positive linear maps
 on $C[0,2\pi]$ of norm less than or equal to $1$.
 We end this section with the proof of two theorems, which are the stronger versions
  of the Theorems (5.3) and (5.4) of \cite{st1}.
  \bt \label{LPO1}
Let $L_n[U_n](g)= g + \epsilon_n(g)$ for every $g$ in the finite set
$\left\{g_1,g_2,g_3,\,.\,.\,.\,.g_m, \sum_{k=1}^m {g_k ^2}\right\}$,
where each $g_k$' s are real valued, continuous
functions and $\epsilon_n(g)$ converges uniformly to 0.
 Then $P_{U_n}(A_n(f))$ converges to $A_n(f)$ in weak cluster for
 all f in the $C^*$- algebra  $\mathbb{A}$ generated by the finite set
 $\left\{g_1,g_2,g_3,\,.\,.\,.\,.g_m\right\}$.
\et
\bp
 First we observe that $L_n[U_n](g)= g + \epsilon_n(g)$ for every $g$ in
  algebra generated by $\left\{g_1,g_2,g_3,\,.\,.\,.\,.g_m\right\},$ by Lemma (\ref{remainder}).
Also we have
\begin{equation} \label{lpo1}
  0 \leq {{\parallel{{A_n}(f_l)- P_{U_n}({A_n}(f_l))}\parallel}_F}^2 = {{\parallel{{A_n}(f_l)}\parallel}_F}^2 -{{\parallel{P_{U_n}{A_n}(f_l)}\parallel}_F}^2
 \end{equation}
 for every function $f_l$ in the algebra generated by $\left\{g_1,g_2,g_3,\,.\,.\,.\,.g_m\right\}$.

 Here ${\parallel (.)\parallel}_F$ denotes the Frobenius norm of matrices.
Also since
 \begin{equation}\nonumber
 L_n[U_n](f_l)= \lambda_i({P_{U_n}({A_n}(f_l))})= f_l({x_i^n }) + \epsilon_n(f_l),
  \end{equation}
 for every $l$, where $\lambda_i( P_{U_n}({A_n}(.)))$ are the eigenvalues of $ P_{U_n}({A_n}(.)),$
 we get the following.
 \begin{equation}\nonumber
 {{\parallel P_{U_n}({A_n}(f_l))\parallel}_F}^2=\sum_{i=1}^{n}{\lambda_i^2({P_{U_n}({A_n}(f_l))})}=
\sum\limits_{i = 1}^n [{\left( {f_l  + \varepsilon _n \left( {f_l } \right)} \right)\left( {x_i^n } \right)]^2 }
 \end{equation}
 Hence
\begin{equation}\nonumber
 {{\parallel P_{U_n}({A_n}(f_l))\parallel}_F}^2=\sum\limits_{i = 1}^n f_l^2\left({x_i^n }  \right)+ o(n).
 \end{equation}
 Since ${\{x_i}^{(n)}\}$ is quasiuniformly distributed (see \cite{st1} for definition), by
Lemma (5.1) in \cite{st1}, we get,
  \begin{equation}\label{lpo2}
  \sum_ {i=0}^{n-1}\left[f_l({x_i}^{(n)}+\varepsilon_n(f_l)({x_i}^{(n)})^2\right]= n/{2\pi}\int_{0}^{2\pi}{f_l^2}+o(n)
  \end{equation}
  Also
  \begin{equation}\nonumber
 {{\parallel {{A_n}(f_l)}\parallel}_F}^2=\sum_{i=1}^{n}{\lambda_i({A_n}(f_l))^2},
 \end{equation}
     for every $l$, and hence by Szego-Tyrtyshnikov Theorem (5.1) in \cite{st1}, we find
  \begin{equation}\label{lpo3}
    {{\parallel {A_n}(f_l)\parallel}_F}^2=n/{2\pi}\int_{0}^{2\pi}{f_l^2} + o(n)
  \end{equation}
  Now from (\ref{lpo1}),(\ref{lpo2}) and (\ref{lpo3}) we get
  \begin{equation} \nonumber
  {{\parallel{{A_n}(f_l)- P_{U_n}({A_n}(f_l))}\parallel}_F}^2 = o(n)
   \end{equation}
 for every function $f_l$ in the algebra generated by $\left\{g_1,g_2,g_3,\,.\,.\,.\,.g_m\right\}$.
   Therefore by Tyrtyshnikov's Lemma (\ref{Tyrtylemma}),
 $P_{U_n}({A_n}(f_l))$ converges to ${A_n}(f_l)$ in weak cluster.
Hence by Theorem (\ref{weakkorovkin}), $P_{U_n}({A_n}(f))$ converges to ${A_n}(f_l)$ in weak cluster
for every f in the $C^*$- algebra  $\mathbb{A}$ generated by
 $\left\{g_1,g_2,g_3,\,.\,.\,.\,.g_m\right\}$.
 Hence the proof is completed.
 \ep
\bt \label{LPO2}
With the assumptions in Theorem (\ref{LPO1}), if $\epsilon_{n}(g) = O(1/n)$ for
 $g$ in the finite set $\{g_1,g_2,\ldots g_m,\sum_{k=1}^m{g_k{^2}}\}$ and
 if the ``grid point algebra'' are uniformly distributed, then the convergence is in strong cluster,
 provided the test functions in the set $\left\{g_1,g_2,g_3,\,.\,.\,.\,.g_m\right\}$
 are Lipschitz continuous and belong to the Krein algebra.
\et
\bp
The proof can be obtained by replacing the polynomials $p$ by
$\left\{g_1,g_2,g_3,\,.\,.\,.\,g_n\right\}$ in the proof of Theorem(5.4) in \cite{st1}.
 The idea is to replace the term of o(n) by constants in the equations (\ref{lpo2}) and (\ref{lpo3}).
 For (\ref{lpo2}), we use the hypothesis $\epsilon_{n}(g) = O(1/n)$ and that
 the ``grid point algebra'' are uniformly distributed.
 For (\ref{lpo3}), we use Widom's theorem (Theorem(5.2) in \cite{st1} or see \cite{Widom89}).
 Hence we will attain
   \begin{equation} \nonumber
  {{\parallel{{A_n}(f_l)- P_{U_n}({A_n}(f_l))}\parallel}_F}^2 = O(1)
   \end{equation}
This completes the proof due to Lemma (\ref{Tyrtylemma}).
\ep
\section{Concluding Remarks:}
We conclude this article by pointing out some important features
of this study. The problems that we encountered in an abstract setting,
includes the link between the spectral information of large matrices, and the
Korovkin-type approximations in the non commutative set up with respect to
various topologies. One can approach the problem in many different dimensions.
 Below, we mention an interesting problem of this regard.
 \begin{itemize}
 \item  The study of  approximation of spectrum of infinite dimensional operators and spectral gap related problems, using truncation method is very  important. But usually these truncations need not be simple. So it would be useful to consider simpler matrices without loosing much of the spectral information. Our aim is to handle this problem using the pre-conditioners and approximating in the sense of clustering of eigenvalues.	
	\end{itemize}

\textbf{Acknowledgments:}
Kiran Kumar is thankful to CSIR, KSCSTE for financial support. Stefano Serra-Capizzano is thankful to the Govt. of Kerala, Erudite program, to the
Italian MiUR, PRIN 2008 N. 20083KLJEZ, for financial support. M.N.N.Namboodiri is thankful to the Govt.Kerala, Erudite program and CUSAT for the support of the realization of this collaborative work.

\end{document}